\documentclass{amsart}

\theoremstyle{plain}
\newtheorem{theorem}{Theorem}
\newtheorem{corollary}[theorem]{Corollary}

\theoremstyle{definition}
\newtheorem{example}{Example}

\begin{document}

\title[Solution of Belousov's problem]{Solution of Belousov's problem}

\author{Maks A. Akivis}
\address{Department of Mathematics,
Jerusalem College of Technology---Mahon Lev,
Havaad Haleumi St., P. O. B. 16031,
 Jerusalem 91160, Israel}
\email{akivis@avoda.jct.ac.il}

\thanks{The research of the
first author was partially supported by the Israel
Ministry of Absorption and the Israel Public Council for Soviet
Jewry.}

\author{Vladislav  V. Goldberg}
\address{ Department of Mathematical Sciences,
New Jersey Institute of Technology,   University Heights, Newark, NJ 07102}
\email{vlgold@m.njit.edu}


\date{March 2000.}

\dedicatory{}



\begin{abstract}
The authors prove that a local $n$-quasigroup defined by the equation
$$
x_{n+1} = F (x_1, \ldots , x_n)
 = \displaystyle \frac{f_1 (x_1) + \ldots + f_n (x_n)}{x_1 + \ldots +
 x_n},
$$
where $f_i (x_i), \; i, j = 1, \ldots , n,$ are
arbitrary functions, is irreducible if and only if
any two functions  $f_i (x_i)$ and  $f_j (x_j), \, i \neq j,$
are not both linear homogeneous, or these functions are linear
homogeneous but   $\frac{f_i (x_i)}{x_i} \neq \frac{f_j (x_j)}{x_j}$.

 This gives a solution of Belousov's problem to construct examples
of irreducible $n$-quasigroups for any $n \geq 3$.

\end{abstract}

\maketitle

\setcounter{equation}{0}

\textbf{1. Introduction} An {\em $n$-quasigroup} is a set $Q$ with
an $n$-ary operation $A$ such that each equation $A (a_1, a_2,
\ldots a_{i - 1}, x, a_{i+1}, \ldots, a_n) = b$ is uniquely
solvable with respect to $x (i = 1, \ldots , n)$.

An $n$-quasigroup $A$, $n > 2$, is {\em reducible} if there exist
 an $k$-ary operation $B$ and $(n-k+1)$-ary
operation $C$ such that $A = B \overset{k}{+} C$.
Otherwise, an $n$-quasigroup $A$ is said to be
{\em irreducible}.

In his monograph \cite{Be}, Belousov  posed the
following problem (p. 217, problem 5):
{\em Construct examples of irreducible $n$-quasigroups, $n > 3$.
Do there exist irreducible $n$-quasigroups for any $n > 3$ ?}
We describe now the development in the solution of this problem.

\begin{itemize}

\item Belousov and Sandik \cite{BS} constructed an example of
an irreducible $3$-quasigroup of order 4 (see also \cite{Be}, p. 115).

\item Frenkin \cite{F} proved  that for any $n \geq 3$,
there exist irreducible $n$-quasigroups of order 4.

\item Using methods of web geometry, Goldberg
\cite{Go1}--\cite{Go2} (see also the book \cite{Go3}, Ch. 4) proved
an existence of infinite local irreducible $n$-quasigroups for any
$n \geq 3$.  It is well-known that the theory of
$(n+1)$-webs is equivalent to the theory
of local differentiable  $n$-quasigroups. Goldberg
proved that in general an arbitrary $(n+1)$-web (or a local
differentiable $n$-quasigroup) is irreducible.

\item One year later, independently, using algebraic methods,
Glukhov \cite{Gl1} (see also \cite{Gl2}) proved an existence of infinite
irreducible $n$-quasigroups for any $n \geq 3$.

\item Borisenko \cite{Bo}  for any $n \geq 3$,
constructed examples of irreducible $n$-quasigroups of finite composite order
$t > 4$.

\end{itemize}

Note that in all these works, {\em no examples of infinite irreducible
 $n$-quasigroups were given.}


 In the current paper we present the {\em simplest} examples
 of local irreducible  and reducible $n$-quasigroups.   They are coordinate $n$-quasigroups
 of a series of irreducible and reducible
 $(n+1)$-webs. We came to these examples from
 the web theory. However, to make this paper accessible for  mathematicians
 not working in web geometry, in our presentation we formulate
 the main results and their proofs without using the web geometry terminology.

 Note also that we could give much more
 examples of local irreducible  $n$-quasigroups but they will be
 more complicated than examples considered in this paper.

\textbf{2. Preliminaries} Suppose that a local differentiable
$n$-quasigroup $A$ is given on a differentiable
manifold $Q$ by the equation
\begin{equation}\label{1}
  x_{n+1} = F (x_1, x_2, \ldots , x_n),
\end{equation}
where $F$ is a $C^2$-function. First, we indicate invariant
conditions for such a local $n$-quasigroup $A$ to be reducible.

Note that for a ternary quasigroup (i.e., when
$n = 3$), Goldberg \cite{Go1} (see also \cite{Go3})
proved that a 3-quasigroup $A$ defined by the equation
\begin{equation}\label{2}
x_4 = F (x_1, x_2, x_3)
\end{equation}
is reducible of type
\begin{equation}\label{3}
x_4 = g (h (x_1, x_2), x_3)
\end{equation}
(where $g$ and $h$ are differentiable functions)
if and only if the function $F$ satisfies the following
second-order nonlinear partial differential equation
\begin{equation}\label{4}
\displaystyle\frac{F_{31}}{F_{32}} = \frac{F_1}{F_2}.
\end{equation}
Here and in what follows, we use the notation
$$
F_i = \frac{\partial F}{\partial x_i}, \;\;
F_{ij} = \frac{\partial^2 F}{\partial x_i \partial x_j}, \;\;
i, j = 1, \ldots , n.
$$

Equation (4)  was noticed by Goursat \cite{Gou}
as far back as   1899 who indicated that function (3), where
$g$ and $h$ are arbitrary functions of two variables each, is
a general solution of equation (4).

We define a local reducible $n$-quasigroup. Without loss of
generality, we can define a local {\em reducible} $n$-quasigroup $A$
as an $n$-quasigroup for which the function $F (x_1, \ldots ,
x_n)$ has the following form:
\begin{equation}\label{5}
x_{n+1} = g (h(x_1, \ldots , x_k), x_{k+1},  \ldots , x_n),
 \end{equation}
(i.e., its operation is reduced to a $k$-ary and $(n + 1 - k)$-ary
operations, $2 \leq k \leq n-1$).
In terminology of \cite{BS} (see also \cite{Be}), such
an $n$-quasigroup is $(1, k)$-reducible. Goldberg (see \cite{Go2}
or \cite{Go3}) found necessary and sufficient
conditions for an $n$-quasigroup (1) to be reducible. For
reducibility of the type (5), these conditions are:
the function $F$ must satisfy the following
system of second-order  nonlinear partial differential equations:
\begin{equation}\label{6}
\displaystyle\frac{F_{pa}}{F_{pb}}
= \frac{F_{a}}{F_b}, \;\; a, b = 1, \ldots k,
\; a \neq b; \; p, q = k+1, \ldots , n.
\end{equation}


The proof of this statement is straightforward: conditions (6)
can be obtained from (5) by differentiation, and it can be shown
that the function defined by equation (5) is a general solution
of the system (6).



Let us consider a few examples.

\begin{example}
If an $n$-quasigroup is (1, 2)-reducible, i.e.,
if
\begin{equation}\label{7}
x_{n+1} = g (h (x_1, x_2), x_3 \ldots , x_n),
\end{equation}
then the conditions (6) take the form
\begin{equation}\label{8}
\displaystyle\frac{F_{p1}}{F_{p2}} = \frac{F_1}{F_2},
\;\;  p = 3, \ldots , n.
\end{equation}
\end{example}

\begin{example}
If an $n$-quasigroup is $(1, 2)$- and $(3, 5)$-reducible, i.e., if
\begin{equation}\label{9}
x_{n+1} = g (h (x_1,  x_2), k (x_3, x_4, x_5), x_6,  \ldots , x_n),
\end{equation}
then the conditions (6) take the form
\begin{equation}\label{10}
\renewcommand{\arraystretch}{1.3}
\left\{
\begin{array}{ll}
\displaystyle\frac{F_{p1}}{F_{p2}} = \frac{F_1}{F_2}, & p = 3, \ldots , n;
  \vspace*{2mm} \\
\displaystyle\frac{F_{\sigma a}}{F_{\sigma b}} = \frac{F_a}{F_b},
&  a, b = 3, 4, 5, \; a \neq b; \; \sigma = 1, 2, 6, 7,  \ldots , n.
\end{array}
\right.
\renewcommand{\arraystretch}{1}
\end{equation}
\end{example}

\begin{example}
If an $n$-quasigroup is $(1, 2, 3)$- and $(1, 2)$-reducible, i.e., if
\begin{equation}\label{11}
x_{n+1} = g (h (k (x_1, x_2), x_3)), x_4,  \ldots , x_n),
\end{equation}
then the conditions (6) take the form
\begin{equation}\label{12}
\renewcommand{\arraystretch}{1.3}
\left\{
\begin{array}{ll}
\displaystyle\frac{F_{pa}}{F_{pb}}= \frac{F_a}{F_b},
\;\; a, b = 1, 2, 3, \; a \neq b; \; p = 4, 5, \ldots , n; \;   \vspace*{2mm}\\
\displaystyle\frac{F_{31}}{F_{32}} = \frac{F_1}{F_2}.
\end{array}
\right.
\renewcommand{\arraystretch}{1}
\end{equation}
\end{example}

These three examples show how to get
 conditions (6) for different kinds
of reducibilities.

Now we will formulate our main theorem.
\begin{theorem}
A local $n$-quasigroup defined by the equation
\begin{equation}\label{13}
x_{n+1} = F (x_1, \ldots , x_n)
 = \displaystyle \frac{f_1 (x_1) + \ldots + f_n (x_n)}{x_1 + \ldots + x_n}
\end{equation}
is irreducible if and only if  any two functions  $f_i (x_i)$
and  $f_j (x_j), \; i \neq j,$
are not both linear homogeneous, or these functions are both linear
homogeneous but   $\frac{f_i (x_i)}{x_i} \neq \frac{f_j (x_j)}{x_j}$.
\end{theorem}

The proof of this theorem follows from the next theorem
on necessary and sufficient conditions for an $n$-quasigroup
defined by equation (13) to be reducible.

Note that the function (13) defines an $n$-cone system
(see \cite{R}).

\begin{theorem}
A local $n$-quasigroup defined by equation $(13)$
is $(1, k)$-reducible if and only if
\begin{equation}\label{14}
f_a (x_a) = c x_a, \;\; a = 1, \ldots , k,
\end{equation}
where $c$ is the same constant for all $a = 1, \ldots , k.$
\end{theorem}

\begin{proof}
A necessary and sufficient condition for an $n$-quasigroup to be
$(1, k)$-reducible is that the function $F$ from (13) satisfies
the system of partial differential equations (6).
Write equations (6) for any two fixed different values $a$ and
$b, \; a, b = 1, \ldots , k$:
\begin{equation}\label{15}
\frac{F_{pa}}{F_{pb}} = \frac{F_a}{F_b}.
\end{equation}
For the $n$-quasigroup defined by equation (13), equation (15)
takes the form
\begin{equation}\label{16}
(f_a' (x_a) - f_b' (x_b)) [f_p' (x_p) (x_1 + \ldots + x_n) -
(f_1 (x_1) + \ldots + f_n (x_n))] = 0.
\end{equation}

First we assume that
$$
f_p' (x_p) (x_1 + \ldots + x_n) - (f_1 (x_1) + \ldots + f_n (x_n))
= 0.
$$
Then
$$
f_p'(x_p) = \frac{f_1 (x_1) + \ldots + f_n (x_n)}{x_1 + \ldots + x_n} (= F).
$$
Since the left-hand side does not depend on $x_a$, then
  $F_a = 0$.
But $F_p = 0$ implies that  $f'_p (x_p) \sum_i x_i - \sum_i f_i (x_i) = 0$,
i.e., $f'_a (x_a) = F$. This along with $f'_p (x_p) = F$  leads to
$f'_i (x_i) = F, \; i = 1, \ldots , n$.
These equalities are possible if and only if
$$
F = A = (\mbox{{\rm const.}})
$$
However, in this case equation (13) does not define
a local $n$-quasigroup since it is not solvable with respect
to the variables $x_i, \; i = 1, \ldots ,n$.

Leaving this case aside, we find from condition
(16) that
$$
f_a' (x_a) = f_b' (x_b).
$$
We can see again that both sides of this equations are constant:
$$
f_a' (x_a) = f_b' (x_b) = c (= \mbox{{\rm const.}}).
$$
Integration gives
$$
f_a (x_a) = c x_a, \;\;  f_b (x_b) = c x_b.
$$
Since $a$ and $b$ are arbitrary numbers from $1, \ldots , k$,
this proves that
$$
f_a (x_a) = c x_a, \;\; a = 1, \ldots , k.
$$

Thus
\begin{equation}\label{17}
F = \frac{c (x_1 + \ldots + x_k) + f_{k+1} (x_{k+1}) +  \ldots + f_n (x_n)}{x_1 +
\ldots + x_n}.
\end{equation}
\end{proof}



The following corollary gives a geometric meaning of local reducible
$n$-quasigroups defined by equation (13). First note that in $\mathbb{R}^n$
equation (13) determines an $(n+1)$-web $W$ formed by  hyperplanes
$x_i = c_i, \; i = 1, \ldots , n$,
parallel to the coordinate hyperplanes of a Cartesian coordinate system
of $\mathbb{R}^n$ and by a family $\lambda_{n+1}$ of hypersurfaces $V$
defined by the equations $F = \alpha (= \mbox{{\rm const.}})$

\begin{corollary}
A local $n$-quasigroup $(13)$ is reducible if and only if in  a Cartesian coordinate
system of  $\mathbb{R}^n$, the normal vector to any hypersurface of the $(n+1)$-web $W$ has
at least two equal coordinates $($i.e., if at least two projections
of this vector onto the coordinate axes  are equal$)$.
\end{corollary}

\begin{proof}
In fact, the equation of a hypersurface $V$ is
$$
f_1 (x_1) + \ldots + f_n (x_n) - \alpha (x_1  + \ldots + x_n) = 0
$$
where $\alpha$ is a constant.

In  a Cartesian coordinate system of  $\mathbb{R}^n$,
the normal vector ${\bf N}$ at an arbitrary point of the
hypersurface $V$ has the coordinates
$$
f_1' (x_1) - \alpha, \ldots , f_n' (x_n) - \alpha.
$$
Thus if two of these coordinates are equal, this implies
$f_i' (x_i) = f_j' (x_j), \; i \neq j$. As we saw in the proof of Theorem 2,
this implies that $f_i (x_i) + f_j (x_j) = c (x_i + x_j)$,
and the local $n$-quasigroup is reducible.

  The converse statement is trivial: it follows from equation
  (17) if one calculate the coordinates of a normal vector.
\end{proof}

By Theorem 2, for reducibilities of types (7), (9),
and (11), the function $F$ defined by (13) has the forms
$$
F = \frac{c (x_1 + x_2) + f_3 (x_3) + \ldots + f_n (x_n)}{x_1 + \ldots +
x_n},
$$
$$
F = \frac{c (x_1 + x_2) + e (x_3 + x_4 + x_5) + f_6 (x_6) +
\ldots + f_n (x_n)}{x_1 + \ldots + x_n},
$$
and
$$
F = \frac{c (x_1 + x_2 + x_3) + f_4 (x_4) +
\ldots + f_n (x_n)}{x_1 + \ldots + x_n}
$$
(where $c$ and $e$ are constants), respectively.

\begin{example}
It follows from Theorem 1 that the simplest example of an irreducible $n$-quasigroup
of type (13) we obtain by taking the  functions $f_i (x_i) = c_i x_i, \; i, j
= 1, \ldots , n$, where $c_i \neq c_j$ if $i \neq j$.

One can produce numerous examples of irreducible quasigroups of this kind.
For example, the  local $n$-quasigroup  defined by the equation
\begin{equation}\label{18}
x_{n+1} = \frac{x_1 +  2 x_2 + 3 x_3 + \ldots + n x_n}{x_1 + \ldots +
x_n}, \;\; n \geq 3,
\end{equation}
\noindent
is irreducible. An $(n+1)$-web in $\mathbb{R}^n$ corresponding to
this $n$-quasigroup is formed by $n$ pencils of hyperplanes
parallel to the coordinate hyperplanes of a Cartesian coordinate
system of  $\mathbb{R}^n$  and a pencil of
hyperplanes whose axis is an $(n-2)$-plane defined by the equations
$$
x_1 +  2 x_2 + 3 x_3 + \ldots + n x_n = 0, \;\;
x_1 + \ldots + x_n = 0.
$$
Note that for $n = 2$, equation (18) defines a parallelizable
web in a 2-plane.
\end{example}

It is easy to see that a local $n$-quasigroup defined by equation (13)
is isotopic to the  $n$-quasigroup defined by the equation
\begin{equation}\label{19}
x_{n+1} = F (x_1, \ldots , x_n)
 = \displaystyle \frac{f_1 (x_1) + \ldots + f_n (x_n) + A}{x_1
 + \ldots + x_n + a},
\end{equation}
where $A$ and $a$ are constants. In fact, we can make two
successive isotopic transformations for the quasigroup (19):
$$
x_n + a \rightarrow x_n.
$$
and
$$
f_n (x_n - a) + A \rightarrow f_n (x_n)
$$
As a result, we will get the $n$-quasigroup defined by
equation (13). Thus the local $n$-quasigroups defined by
equations (13) and (19) are isotopic.

\begin{example} If we take $f (x_i) = x_i^2$ in (21), then
we obtain the local $n$-quasigroup defined by
\begin{equation}\label{20}
x_{n+1} = F (x_1, \ldots , x_n)
 = \displaystyle \frac{(x_1)^2 + \ldots + (x_n)^2 + A}{x_1
 + \ldots + x_n + a}.
\end{equation}
By Theorem 1, this  $n$-quasigroup is irreducible.

It is easy to see that $F = \alpha (= \mbox{{\rm const}})$
defines a family $\lambda_{n+1}$ of
hyperspheres in $\mathbb{R}^n$. If all these hyperspheres pass through the points
$(1, 0, \ldots , 0), (0, 1, \ldots , 0), \linebreak (0, 0, \ldots , 1)$, then
(20) implies that $1 + A = \alpha (1 + a)$. Since the last
equation must be valid for any $\alpha$,
it follows  that $A = a = -1$. The family $\lambda_{n + 1}$
along with the families $\lambda_i, \; i = 1, \ldots , n,$
of  hyperplanes $x_i = c_i$   parallel to the coordinate
 hyperplanes of a Cartesian coordinate system of
 $\mathbb{R}^n$ form an irreducible $(n+1)$-web.

In particular, if $n = 3$, we have
$$
x_4 = \displaystyle \frac{(x_1)^2 + (x_2)^2 + (x_3)^2 - 1}{x_1
 + x_2 + x_3 -1}.
$$
In this case the family $\lambda_4$ of the web $W$
is the  1-parameter family of hyperspheres in $\mathbb{R}^4$
passing through the points (1, 0, 0), (0, 1, 0), and (0, 0, 1).

Moreover if $n = 2$, then we have
$$
x_3 = \displaystyle \frac{(x_1)^2 + (x_2)^2 - 1}{x_1 + x_2 -1}.
$$
In this case the equation
$$
(x_1)^2 + (x_2)^2 - 1 - \alpha (x_1 + x_2 -1) = 0, \; \alpha =
\mbox{{\rm const.}}
$$
 determines an 1-parameter family of circles in $\mathbb{R}^3$
passing through the points (1, 0) and (0, 1).
This family and two 1-parameter families of straight lines
parallel to the coordinate lines of a Cartesian coordinate system of
 $\mathbb{R}^2$ form a nonhexagonal 3-web (cf. \cite{Bl},
\S 3, where the same  3-web is considered---the only difference is
that in \cite{Bl}, circles pass through the points (0, 0) and (1, 1)).
\end{example}

\makeatletter \renewcommand{\@biblabel}[1]{\hfill#1.}\makeatother
\bibliographystyle{amsplain}

\end{document}